\documentstyle[12pt]{article} 
\begin{document}
\newcommand{\singlespace}{
    \renewcommand{\baselinestretch}{1}
\large\normalsize}
\newcommand{\doublespace}{
   \renewcommand{\baselinestretch}{1.2}
   \large\normalsize}
\renewcommand{\theequation}{\thesection.\arabic{equation}}

\input amssym.def
\input amssym
\setcounter{equation}{0}
\def \ten#1{_{{}_{\scriptstyle#1}}}
\def \Z{\Bbb Z}
\def \C{\Bbb C}
\def \R{\Bbb R}
\def \Q{\Bbb Q}
\def \N{\Bbb N}
\def \l{\lambda}
\def \V{V^{\natural}}
\def \wt{{\rm wt}}
\def \tr{{\rm tr}}
\def \Res{{\rm Res}}
\def \End{{\rm End}}
\def \Aut{{\rm Aut}}
\def \mod{{\rm mod}}
\def \Hom{{\rm Hom}}
\def \im{{\rm im}}
\def \<{\langle} 
\def \>{\rangle} 
\def \w{\omega}
\def \c{{\tilde{c}}}
\def \o{\omega}
\def \t{\tau }
\def \ch{{\rm ch}}
\def \a{\alpha }
\def \b{\beta}
\def \e{\epsilon }
\def \la{\lambda }
\def \om{\omega }
\def \O{\Omega}
\def \qed{\mbox{ $\square$}}
\def \pf{\noindent {\bf Proof: \,}}
\def \voa{vertex operator algebra\ }
\def \voas{vertex operator algebras\ }
\def \p{\partial}
\def \1{{\bf 1}}
\def \ll{{\tilde{\lambda}}}
\def \H{{\bf H}}
\def \F{{\bf F}}
\def \h{{\frak h}}
\def \g{{\frak g}}
\def \rank{{\rm rank}}
\singlespace
\newtheorem{thmm}{Theorem}
\newtheorem{co}[thmm]{Corollary}
\newtheorem{th}{Theorem}[section]
\newtheorem{prop}[th]{Proposition}
\newtheorem{coro}[th]{Corollary}
\newtheorem{lem}[th]{Lemma}
\newtheorem{rem}[th]{Remark}
\newtheorem{de}[th]{Definition}
\newtheorem{con}[th]{Conjecture}
\newtheorem{ex}[th]{Example}

\begin{center}
{\Large {\bf  Toward classfication of rational vertex operator algebras 
with central charges less than 1}} \\
\vspace{0.5cm}

Chongying Dong\footnote{Supported by NSF grants and a research grant from the
Committee on Research, UC Santa Cruz}  \  and \ 
Wei Zhang\\
Department of mathematics\\
University of California\\
 Santa Cruz, CA 95064

\end{center}
\hspace{1.5 cm}

\begin{abstract} 
The rational and $C_2$-cofinite simple vertex operator algebras whose
effective central charges $\title{c}$ and the central charges $c$ are
equal and less than 1 are classified. Such a vertex operator algebra is
zero if $\c<0$ and $\C$ if $\c=0.$ If $\c>0,$ it is an extension of
discrete Virasoro vertex operator algebra $L(c_{p,q},0)$ by its
irreducible modules. It is also proved that for any rational and 
$C_2$-cofinite simple vertex operator algebra with $c=\c,$ the vertex operator
subalgebra generated by the Virasoro vector is simple.
\end{abstract}

\section{Introduction}

One of the most important problems in the theory of vertex operator
algebra is to classify the rational vertex operator algebras. It is
not realistic to achieve this goal at this stage due to the limited
knowledge of the structure theory and representation theory. 
If a vertex operator algebra is rational then the
central charge $c$ and effective central charge $\tilde{c}$ are
rational (cf. \cite{AM}, \cite{DLM2}). While the central
charge can be negative, the effective central charge is always
nonnegative \cite{DM2}. In this paper we classify the rational 
vertex operator algebras with $c=\tilde{c}<1$ although we cannot write down
the results explicitly.

It is well known that one can construct vertex operator algebras
associated to highest weight modules for the Virasoro algebra \cite{FZ}.
In particular, each irreducible highest weight module $L(c,0)$ for any complex
number $c$ is a simple vertex operator algebra. Moreover, $L(c,0)$ is rational
if and only if $c=c_{p,q}=1-6(p-q)^2/pq$ for coprime positive integers
$p,q$ with $1<p,q$ \cite{W}. Furthermore, $L(c,0)$ is unitary if and only if
$c=c_{p,p+1}$ for all $p>1$ or $c\geq 1$ (see \cite{FQS} and \cite{GKO}).
Our classification result says that any simple, rational and $C_2$-cofinite
vertex operator algebra with $c=\c<1$ is an extension of
the Virasoro vertex operator algebra $L(c_{p,q},0)$ for some $p,q.$ 
That is, such vertex operator algebra is a finite direct sum of irreducible  
$L(c_{p,q},0)$-modules. 
 
The main ideal is to use  the modular invariance of the $q$-characters of the 
modules (see \cite{Z} and \cite{DLM2})  to control the growth of the graded
dimensions of the vertex operator algebra. The same idea has been used to 
classify the holomorphic vertex
operator algebras with small central charges \cite{DM1}, to prove the
nonnegative property of the effective central charges \cite{DM2},
and to obtain some uniqueness result on the moonshine vertex operator algebra
$V^{\natural}$ \cite{DGL}. The modular invariance property of 
the $q$-characters of the modules is also the reason why we use the effective 
central charges instead of central charges (see Lemma \ref{l2.1} below).
We should point out that we do not assume that the vertex operator
algebra is a unitary representation for the Virasoro algebra and that 
$\sum_{i}|\chi_i(q)|^2$ is a modular function over the full modular group
where $\chi_i(q)$ are the $q$-characters of the irreducible modules of
the vertex operator algebra (see Section 2).  

As a corollary of the main result we prove that for any simple, rational
and $C_2$-cofinite vertex operator algebra with $c=\c$, the vertex operator subalgebra 
generated by the Virasoro vector is simple. That is, the vertex operator 
subalgebra is an irreducible  highest weight module for the Virasoro algebra.

It is worthy to mention that we do not have a explicit list of such
vertex operator algebras. An eventual classification
requires to construct all extensions of $L(c_{p,q},0)$ for all $p,q.$
In the case that
$c=c_{p,p+1}$, the extensions of $L(c_{p,p+1},0)$ have been classified
in the theory of conformal nets (an analytical approach to conformal
field theory) \cite{KL} (also see \cite{X}). Although it is believed that such
classification result is valid in the theory of vertex operator
algebra, most of such extensions have not been constructed in the context of
vertex operator algebra except for a few examples from the code vertex operator
algebras and lattice vertex operator algebras.

\section{Rational vertex operator algebras}
\setcounter{equation}{0}

In this section, we review some basic facts on the $q$-characters 
of modules for a rational vertex operator algebra. The main
feature of these functions is the modular invariance property
\cite{Z}, and its connection with the vector-valued modular forms
\cite{KM}. This connection is the key for us to estimate the growth of
the graded dimensions of the vertex operator algebra and its modules.
We will also discuss the effective central charge $\c.$ 

We assume that vertex operator algebra $V$ is simple and
is of CFT type. That is,
\begin{equation}\label{2.1}
V = \bigoplus_{n=0}^{\infty}V_n
\end{equation} 
moreover $V_0$ is spanned by the vacuum vector $\1$. Following \cite{DLM1}, $V$ is called rational if the admissible module category is semisimple. $V$ is
called $C_2$-cofinite if $V/C_2(V)$ is finite dimensional \cite{Z}
where $C_2(V)=\<u_{-2}v|u,v\in V\>.$

A rational  vertex operator algebra $V$ has only 
finitely many irreducible modules  $V = M^{1}, M^{2}, ... M^{r}$ up to 
isomorphism such that 
$$M^i=\oplus_{n\geq 0}M^i_{\l_i+n}$$
where $\l_i$ is a rational number and $M^i_{\l_i}\ne 0$ (see \cite{DLM1},
\cite{DLM2}). Moreover each homogeneous subspace 
$M^i_{\l_i+n}$ is finite dimensional. 
Let $\l_{min}$ be the minumum among the $\l_i.$ The
 effective central charge $\tilde{c}$ which appeared 
in the physics literature [GN] is defined by $\c=c - 24 \lambda_{min}.$ One of the main results in \cite{DM2} is that
$\tilde{c}$ is nonnegative, and $\tilde{c}=0$ if and only if $V=\C$ is trivial.

For each $i$ we define the $q$-character of $M^i$ as 
$$\chi_i(q)=ch_qM^i=\tr_{M^i}q^{L(0)-c/24}=\sum_{n\geq 0}(\dim M^i_{n+\l_i})q^{n-c/24}.$$ It is proved in \cite{Z} (also see \cite{DLM2}) that if $V$ is
rational and $C_2$-cofinite then each $\chi_i(q)$ is
a holomorphic function on the upper half plane ${\Bbb H}$
where $q=e^{2\pi i\tau}$ and the span
of these functions affords a representation of the modular group $SL(2,\Z).$
For short we also write $\chi_i(\tau)$ for $\chi_i(q).$ 
Then there exists a group homomorphism $\rho$ from $SL(2,\Z)$ to $GL(r,\C)$ 
such that for any $\gamma\in SL(2,\Z),$ 
$$\chi_i(\gamma\tau)=\sum_j\gamma_{ij}\chi_j(\tau)$$
where $\rho(\gamma)=(\gamma_{ij}).$ 
This exactly says that  $\chi(\tau)=(\chi_1(\tau),\cdots,\chi_r(\tau))$ 
is a ( meromorphic) vector-valued modular function \cite{KM}. 

Recall the Dedekind eta function
$$\eta(\tau) = q^{1/24}\phi(q) = q^{1/24}\prod_{n=1}^{\infty}(1 - q^n)$$
 and the expansion
$$\frac{1}{\prod_{n\geq 1}(1-q^n)} =\sum_{n=0}^{\infty}p(n) q^n$$
where $p(n)$ is the usual unrestricted partition function. An asymptotic
expression for $p(n)$ is given by 
$$p(n)\sim \frac{e^{\pi\sqrt{2n/3}}}{4n\sqrt{3}}$$
as $n\to \infty$ (cf. \cite{A}).
It is clear that $p(n)$ 
grows faster than  $n^{\alpha}$ for any fixed real number $\alpha.$

The $\eta(\tau)$ is a modular form of weight $1/2.$ Since
$\eta(\tau)^{\tilde{c}}\chi_i(\tau)$ is 
holomorphic at $\tau=i\infty,$  
 $$\eta(\tau)^{\tilde{c}}\chi(\tau)=(\eta(\tau)^{\tilde{c}}\chi_1(\tau),\cdots,\eta(\tau)^{\tilde{c}}\chi_r(\tau))$$ 
is a holomorphic  vector-valued modular form of weight $\tilde{c}/2.$  
From \cite{KM} the Fourier 
coefficients $a_n$ of a component of a 
holomorphic vector-valued modular form
satisfy the growth condition $a_n = O(n^{\alpha})$ 
for a constant $\alpha$ independent of $n$. As a result, we see that
\begin{lem}\label{l2.1} The Fourier coefficients of each component of 
$\eta(\tau)^{\c}\chi_(\tau)$ satisfy a polynomial growth condition 
$a_n = O(n^{\alpha})$.
\end{lem}

\section{Vertex operator algebras with $c<1$}
\setcounter{equation}{0}

In this section we will prove the main result. We assume that
$V$ is a rational and $C_2$-cofinite vertex operator algebra with 
$c=\tilde{c}<1.$ It is proved in \cite{DM2} that $\c$ is always 
nonnegative. If $c=0$ then $V=\C.$ So we assume that $c>0.$ 

We first need information on the vertex
operator algebras associated to the highest weight modules for the Virasoro
algebra (see \cite{FF}, \cite{FZ}, \cite{W}). 
We use the standard basis $\{L_n,C|n\in\Z\}$ for the Virasoro algebra. 
For any two complex numbers $c,h$ we denote the Verma module with
central charge $c$ and highest weight $h$ by $V(c,h),$ as usual. 
Let $\bar{V}(c,0)$ be the quotient of $V(c,0)$ modulo submodule 
generated by $L_{-1}v$ where $v$ a nonzero highest weight vector 
of $V(c,0)$ with highest weight $0.$ We use $L(c,h)$ to denote the irreducible
quotient of $V(c,h).$ 

We have already defined the $q$-character $\ch_qM$ if $M$ is a module for any
vertex operator algebra. We now extend this definition to any module
for the Virasoro algebra with finite dimensional homogeneous subspaces
and central charge $c.$ In general the $q$-character is just a formal power
series in $q.$ Note that $\ch_q \bar{V}(c, 0)=\frac{q^{-c/24}}{\prod_{n>1}(1-q^n)}$ and its coefficients grow faster than $n^{\alpha}$ for any fixed
real number $\alpha.$

\begin{lem}\label{la} For any
$\mu>0$ the coefficients of $\frac{1}{\prod_{n>1}(1-q^n)^{\mu}}$ grow faster
than any polynomial $n^\alpha.$ 
\end{lem}

\pf Observe that the coefficients ${-\mu\choose i}(-1)^i$ of $q^i$ in  the expansion of $(1-q)^{-\mu}$ is always positive for any $\mu >0.$ Assume that the coefficients of 
$$\frac{1}{\prod_{n>1}(1-q^n)^{\mu}}=\sum_{n\geq 0}a_nq^n$$
satisfy the polynomial growth condition. Then 
the coefficients of $\frac{1}{\prod_{n>1}(1-q^n)^{k\mu}}$ satisfy the 
polynomial growth condition for any positive integer $k.$ But if $k$ is large 
enough then $k\mu>1$ and the coefficients of
 $\frac{1}{\prod_{n>1}(1-q^n)^{k\mu}}$ grow faster than
$n^{\alpha}$ for any real number $\alpha.$ 
 This is a contradiction. \qed
\bigskip

Here are some basic facts about these modules.
\cite{FF}, \cite{FQS},\cite{GKO}, \cite{FZ}, \cite{W}.

\begin{prop}\label{vir} Let $c$ be a complex number. Then the following hold:

(i) $\bar{V}(c,0)$ is a vertex operator algebra and $L(c,0)$ is a simple
vertex operator algebra.

(ii) The following are equivalent: (a) $\bar{V}(c, 0) = L(c, 0),$ 
(b) $c\ne c_{p,q}=1-6(p-q)^2/pq$ for all coprime 
positive integers $p,q$ with $1<p<q,$ 
(c)  $L(c,0)$ is not rational. 
In this case, the $q$-graded character of $L(c,0)$ is equal 
to $\frac{q^{-c/24}}{\prod_{n>1}(1-q^n)}$ and its
coefficients grow faster than any polynomial in $n.$ 

(iii) The following are equivalent: (a) $\bar{V}(c, 0) \ne L(c, 0),$ 
(b) $c=c_{p,q}$ for some $p,q,$ (c) $L(c,0)$ is rational.  
\end{prop}

We now back to vertex operator algebra $V.$ 
Let $U=\<\omega\>$ be the vertex operator subalgebra of $V.$ Then there are
two possibilities: either $U$ is isomorphic to $\bar{V}(c,0)$ or $L(c,0)$
from the structure theory for these modules \cite{FF}.

The following is the key lemma.
\begin{lem}\label{kl} Assume that $c<1.$ Then the $q$ character $ch_qU$ of $U$ is different 
from $\frac{q^{-c/24}}{\prod_{n>1}(1-q^n)}.$ 
\end{lem}

\pf We prove by contradiction. Suppose that
$\ch_qU$ is equal to $\frac{q^{-c/24}}{\prod_{n>1}(1-q^n)}.$ Then
$$ \eta(q)^c \ch_qU =\frac{(1-q)^c}{\prod_{n>1}(1-q^n)^{1-c}}.$$
By Lemma \ref{la} the coefficients of $\frac{1}{\prod_{n>1}(1-q^n)^{1-c}}$ grow faster than
any polynomial in $n.$ Set
$$\eta(q)^c\ch_qV=\sum_{n\geq 0}b_nq^n.$$ By Lemma \ref{l2.1}, the
coefficients $b_n$ satisfy polynomial growth condition.

 Let
$f(q)=(1-q)^{-1}\eta(q)^c\ch_qV=\sum_{n\geq 0}c_nq^n.$ Then
$$c_n=\sum_{i=0}^nb_i$$ for all $n\geq 0.$ Since $b_n$ satisfy 
polynomial growth condition, there exist positive constants $C$ and $\alpha$
such that $|b_n|\leq Cn^{\alpha}$ for all $n.$ As a result,
$|c_n|\leq C(n+1)n^{\alpha}\leq 2Cn^{\alpha+1}$ for $n>0.$ 
That is, the coefficients
$c_n$ also satisfy the polynomial growth condition.
 
Since $U$ is a subspace of $V$ we see that $\dim U_n\leq \dim V_n$ for all
$n.$ This implies that   $\ch_qU\leq \ch_qV$ and 
$\eta(q)^c \ch_qU\leq \eta(q)^c\ch_qV$ as real numbers 
for any $q\in (0,1).$  Thus $\frac{(1-q)^c}{\prod_{n>0}(1-q^n)^{1-c}}\leq f(q)$
for all $q\in (0,1).$ But this is impossible as the coefficients
of $\frac{1}{\prod_{n>1}(1-q^n)^{1-c}}$ grow faster than any polynomial 
in $n.$
The proof is complete. \qed

\begin{coro}\label{c} Let $V$ be a simple, rational and $C_2$-cofinite 
vertex operator
algebra with central charge $c=\tilde{c}<1.$ Then the vertex operator
algebra $U$ generated by the Virasoro element $\omega$ is simple and
$c=c_{p,q}$ for some coprime $p,q$ such that $1<p<q.$
\end{coro}

\pf  If $c\ne c_{p,q}$ then $\bar{V}(c,0)$ is irreducible module for the 
Virasoro algebra by Proposition \ref{vir}. This implies that $U$ is isomorphic to
$\bar{V}(c,0)$ and 
$$ch_qU=\ch_q\bar{V}(c,0)=\frac{q^{-c/24}}{\prod_{n>1}(1-q^n)}$$
which contradicts to Lemma \ref{kl}. 

So we can assume that $c=c_{p,q}$ for coprime $p,q$ such that $1<p<q.$
In this case $U$ is either isomorphic to $\bar{V}(c_{p,q},0)$ or
$L(c_{p,q},0)$ as a module for the Virasoro algebra. Again by Lemma
\ref{kl}, $U$ is not isomorphic to $\bar{V}(c_{p,q},0)$ by grading
restriction. This forces $U$ to be isomorphic to $L(c_{p,q},0),$ as
desired.  \qed
\bigskip

A vertex operator algebra $V^1$ is called an {\em extension}
of another vertex operator algebra $V^2$ if $V^2$ is isomorphic to a vertex
operator subalgebra of $V^1$ with the same Virasoro vector. In particular,
$V^1$ is a $V^2$-module. Note that 
if $V^2$ is rational then $V^2$ has only finitely many irreducible modules
up to isomorphism and $V^1$ is necessarily a finite direct sum of irreducible
$V^2$-module due to finiteness of the homogeneous subspaces of vertex
operator algebras.

We are now in a position to present the main result of this paper.
\begin{th}\label{t1} Let $V$ be a simple, rational and $C_2$-cofinite 
vertex operator
algebra with central charge $c=\tilde{c}<1.$ Then $V$ is an extension
of $L(c_{p,q},0)$ by its simple modules for some coprime $p,q$ such
that $1<p<q.$
\end{th}

\pf By Proposition \ref{vir}, the vertex operator subalgebra
$U$ of $V$ is rational. So as module for $U$ $V$ is a direct sum of 
finitely many irreducible $U$-modules. \qed
\bigskip

So this theorem reduces the classification of simple, rational and
$C_2$-cofinite vertex operator algebras with central charge
$c=\tilde{c}<1$ to the problem of classification of extensions of
vertex operator algebras $L(c_{p,q},0).$ 
It  will certainly be a difficult and complicated task to
construct the extensions of $L(c_{p,q},0)$ for all $p,q.$ One
definitely needs both coset construction and orbifold theory to
achieve this. We will not go to this direction in this paper.
See \cite{LLY} for some extensions of $L(c_{p,p+1},0)$ by a simple
module.

We remark that the condition that $c=\c$ in Theorem \ref{t1} is
necessary.  Here is a counter example in which 
$V$ is a simple, rational
and $C_2$-cofinite vertex operator algebra with $c<1,$ but $V$ is not an extension
of $L(c_{p,q},0)$ for any coprime integers $1<p<q.$ Recall that $L(c_{2,5},0)$ is a
simple, rational and $C_2$ cofinite vertex operator algebra with
central charge $-22/5$ and effective central charge $2/5.$ Let $W$ be
any simple, rational and $C_2$-cofinite vertex operator algebra whose
central charge and effective central charge are $5.$ For example, one
can take $W$ to be a lattice vertex operator algebra $V_L$ for a rank
5 positive definite even lattice $L.$ Set $V=L(c_{2,5},0)\otimes W.$
Then $V$ is a simple, rational and $C_2$-cofinite vertex operator
algebra with central charge $c=3/5.$ Since any irreducible
$V$-module is a tensor product of an irreducible $L(c_{2,5},0)$-module
with an irreducible $W$-module (cf. \cite{FHL}, \cite{DMZ}), we see
that the effective central charge of $V$ is $\c=27/5.$ Then $3/5$ is
different from $c_{p,q}$ for all coprime integers $p,q$ such that
$1<p<q.$ In order to see this let $c_{p,q}=3/5.$ Then $15(p-q)^2=pq.$
Clearly, $q-p\ne 1.$ Let $d$ be a prime divisor of $q-p.$ Then $d$ is
a also a divisor of $p$ or $q.$ As a result $d$ is a common divisor of
$p$ and $q.$ This is a contradiction. This shows that $V$ is not an
extension of $L(c_{p,q},0)$ by finitely many irreducible modules for
any such $p,q.$

Recall from \cite{DM2} that for vertex operator algebra $L(c_{p,q},0),$
the effective central charge $\c=1-\frac{6}{pq}.$ So 
$c=c_{p,q}=1-\frac{6(p-q)^2}{pq}=\c$ if and only if
$q-p=1,$ or if and only if $L(c_{p,q},0)$ is unitary. So our theorem
eliminate the nonunitary rational Virasoro vertex operator algebras
$L(c_{p,q},0)$ for $q\ne p+1$ in our consideration. But one can still ask
if any extension of $L(c_{p,q},0)$ for $q\ne p+1$ satisfies the assumptions
in Theorem \ref{t1}. Note that that all the irreducible modules
of $L(c_{p,q},0)$ are $L(c_{p,q},h_{r,s})$ where 
$$h_{r,s}=\frac{(sp-rq)^2-(p-q)^2}{4pq}$$
for integers $r,s$ in the ranges $1\leq r\leq p-1$ and $1\leq s\leq q-1.$
So each extension $V$ of $L(c_{p,q},0)$ is a finite direct sum of irreducible
$L(c_{p,q},0)$-module and any irreducible $V$-module is also a finite direct 
sum of   $L(c_{p,q},h_{r,s})'$s. We suspect that the effective
central charge $\c=c$ for $V$ if and only if $c=c_{p,p+1}.$ But we do not have 
a proof for this.

We now turn our attention to an arbitrary simple vertex operator
algebra $V.$ We still denote by $U$ the vertex operator subalgebra
generated by the Virasoro vector $\omega.$ It has been a problem in
the theory of vertex operator algebra whether $U$ is a simple vertex
operator algebra of $U.$ Here is a solution.

\begin{th}\label{t2} If $V$ is a simple, rational and $C_2$-cofinite vertex 
operator algebra. Then the Virasoro vertex operator algebra
$U$ is simple if either $c\geq 1$ or $c=\c<1.$
\end{th}

\pf If $c\ne c_{p,q}$ for all $p,q,$ $\bar{V}(c,0)$ is simple (see
Proposition \ref{vir} and must be isomorphic to $U.$ If $c=c_{p,q},$ $U$ is 
isomorphic to $L(c_{p,q},0)$ by Corollary \ref{c}. \qed
\bigskip

We again remark that if $c\ne \c,$ Theorem \ref{t2} is not valid. In
this case we consider $V=L(c_{2,5},0)^{\otimes 5}\otimes W$ where $W$
is any simple, rational and $C_2$-cofinite vertex operator algebra
such that both the central charge and the effective central charge are
$22.$ Then $V$ is a simple, rational and $C_2$-cofinite vertex
operator algebra with central charge $0$ and effective central charge
$2.$ Let $\omega^1$ and $\omega^2$ be the Virasoro vectors of
$L(c_{2,5},0)^{\otimes 5}$ and $W,$ repectively. Then the Virasoro
vector of $V$ is $\omega=\omega^1\otimes 1+1\otimes \omega^2.$ Clearly,
$\omega\ne 0.$ This implies that the vertex operator subalgebra $U$
genereated by $\omega$ is not equal to $\C=L(0,0).$ So $U$ is not
simple.


\begin{thebibliography}{ABCDE} 
\bibitem[AM]{AM} G. Anderson and G. Moore, Rationality in conformal
field theory, {\em Comm. Math. Phys.} {\bf 117} (1988), 441-450.

\bibitem[A]{A} Tom Apostol, Modular Functions and Direchlet Series in Number Theory, Springer-Verlag, New York, 1990.

\bibitem[DGL]{DGL}
C. Dong, R. Griess Jr. and C. Lam, Uniqueness results of the moonshine
vertex operator algebra, {\em American Journal of Math.} {\bf 129}
(2007), 583-609.

\bibitem[DLM1]{DLM1} C. Dong, H. Li and G. Mason,  Twisted representations of
vertex operator algebras, {\em Math. Ann.} {\bf 310} (1998), 571-600.

\bibitem[DLM2]{DLM2} C. Dong, H. Li and G. Mason, Modular invariance of trace functions in orbifold theory and generalized 
moonshine, {\em Comm. Math. Phys.} {\bf 214} (2000), 1-56.

\bibitem[DLMM]{DLMMM} C. Dong, H. Li, G. Mason and P. Montague,
The radical of a vertex operator algebra,
in: {\em Proc. of the Conference on the Monster and Lie algebras at 
The Ohio State University, May 1996,}  ed. by J. Ferrar and K. Harada,
Walter de Gruyter, Berlin-
New York, 1998, 17-25.
York, 1998.

\bibitem[DM1]{DM1} C. Dong and G. Mason, Holomorphic vertex operator
algebras of small central charges, {\em Pacific J. Math.} {\bf 213} (2004), 
253-266. 

\bibitem[DM2]{DM2} C. Dong and G. Mason,  Rational vertex operator algebras and the
effective central charge, {\em International Math. Research
Notices} {\bf 56} (2004),  2989-3008.

\bibitem[DMZ]{DMZ} C. Dong, G. Mason and Y. Zhu,  Discrete series of the 
Virasoro algebra and the moonshine module, {\em Proc. Symp. Pure. Math., American Math. Soc.} {\bf 56} II (1994), 295-316.


\bibitem[FF]{FF} B. Feigin and D. Fuchs, Verma Modules over the Virasoro Algebra, 
{\em Lect. Notes in Math.} {\bf 1060}, Springer, (1984),  230-245.

\bibitem[FHL]{FHL} I. B. Frenkel, Y. Huang and J. Lepowsky, On
axiomatic approaches to vertex operator algebras and modules,
{\it Memoirs American Math. Soc.} {\bf 104}, 1993.

\bibitem[FLM]{FLM} I. B. Frenkel, J. Lepowsky and A. Meurman, 
Vertex Operator Algebras and the Monster, {\em Pure and Applied
Math.,} Vol. {\bf 134}, Academic Press, 1988.

\bibitem[FZ]{FZ} I. B. Frenkel and Y. Zhu, Vertex operator algebras 
associated to  representations of affine and Virasoro algebras, {\em Duke
 Math. J.} {\bf 66} (1992), 123-168.

\bibitem[FQS]{FQS} D. Friedan, Z. Qiu and S. Shenker, 
Details of the non-unitarity proof for highest weight representations of
Virasoro Algebra, {\em  Comm. Math. Phys.} {\bf 107} (1986), 535-542.

\bibitem[GN]{GN} M. Gaberdiel and A. Neitzke, Rationality, quasirationality
and finite W-algebras, {\em Comm. Math. Phys.} {\bf 238} (2003), 305-331.

\bibitem[GKO]{GKO} P. Goddard, A. Kent and D. Olive,
Unitary representations of the Virasoro Algebra and super-Virasoro algebras,
{\em Comm. Math. Phys.} {\bf 103} (1986), 105-119.

\bibitem[KL]{KL} Y.  Kawahigashi and R. Longo, Classification of local 
conformal nets. Case $c<1,$ {\em  Ann. of Math.} {\bf 160} (2004), 493--522. 

\bibitem[KM]{KM} M. Knopp and G. Mason, On vector-valued modular forms and their Fourier coefficients, {\em Acta Arith.} {\bf 110} (2003), 117-124.

\bibitem[LLY]{LLY} C. Lam, N. Lam and H. Yamauchi, Extensions of unitary
Virasoro vertex operator algebra by a simple module, {\em International Math. Research
Notices} {\bf 2003} (2003), 577-611.


\bibitem[W]{W} W. Wang, Rationality of Virasoro vertex operator algebras, 
{\it International Math. Research Notices}, {\bf 71} (1993), 197-211.

\bibitem[X]{X} F. Xu, New braided endmorphisms from conformal inclusions,
{\em Comm. Math. Phys.} {\bf 192} (1998), 347-403.

\bibitem[Z]{Z} Y. Zhu, Modular invariance of characters of vertex operator algebras,
{\em J. Amer, Math. Soc.} {\bf 9} (1996), 237-302.
\end{thebibliography}
\end{document}